\documentclass[11pt,a4paper,reqno]{amsart} 

\usepackage{amssymb,graphicx}

\usepackage{amssymb}

\newcommand{\koniec}{\begin{flushright}  $\Box $ \end{flushright}}
\newtheorem{theo}{Theorem}[section] 
\newtheorem{prop}[theo]{Proposition}  
\newtheorem{lemma}[theo]{Lemma}

\newcounter{mnotecount}[section]

\renewcommand{\themnotecount}{\thesection.\arabic{mnotecount}}

\newcommand{\mnote}[1]
{\protect{\stepcounter{mnotecount}}$^{\mbox{\footnotesize
$
\bullet$\themnotecount}}$ \marginpar{
\raggedright\tiny\em
$\!\!\!\!\!\!\,\bullet$\themnotecount: #1} }

\newcommand{\CP}{\mathbb{CP}}
\newcommand{\C}{\mathbb{C}}

\newcommand{\R}{\mathbb{R}}
\newcommand{\Rho}{\mathrm{P}}
\def\p{\partial}
\def\be{\begin{equation}}
\def\ee{\end{equation}}
\def\ep{{\varepsilon}}
\def\bea{\begin{eqnarray}}
\def\eea{\end{eqnarray}}

\def\hw{\hat\omega}
\def\hnabla{\hat\nabla}
\def\hk{\hat K}
\def\hrho{\hat\rho}

\def\hS{\hat\Sigma}
\newcommand{\spp}{\mathbb{S}}
\let\a=\alpha

\def\ov{\overline}

\begin{document}
\title{FOUR--DIMENSIONAL METRICS CONFORMAL TO  K\"AHLER}
\author{Maciej Dunajski}
\address{Department of Applied Mathematics and Theoretical Physics\\ 
University of Cambridge\\ Wilberforce Road, Cambridge CB3 0WA\\ UK.}
\email{m.dunajski@damtp.cam.ac.uk}
\author{Paul Tod}
\address{The Mathematical Institute\\
Oxford University\\
24-29 St Giles, Oxford OX1 3LB\\ UK.}
\email{paul.tod@sjc.ox.ac.uk}

\begin{abstract} 
We  derive some necessary conditions on a Riemannian metric $(M, g)$ in 
four dimensions 
for it to be locally conformal to K\"ahler. If the conformal curvature is
non anti--self--dual,  the self--dual  Weyl spinor must be 
of algebraic type 
$D$ and satisfy a simple first order conformally invariant 
condition which is necessary and
sufficient   for the existence of a K\"ahler metric in the conformal class.
In the anti--self--dual case we establish 
a one to one correspondence between K\"ahler
metrics in the conformal class and non--zero 
parallel sections of a certain connection on a natural rank ten
vector bundle over $M$.
We use this characterisation to  provide examples of ASD metrics 
which are not conformal to K\"ahler.

We establish a  link between
the `conformal to K\"ahler condition' in dimension four
and the metrisability of projective structures in dimension two.
A projective structure on a surface $U$ is metrisable if and only if the  induced (2, 2) conformal structure on $M=TU$ admits a K\"ahler metric or a para--K\"ahler metric.

\end{abstract}   
\maketitle
\section{Introduction}
A K\"ahler structure on a $2n$--dimensional real manifold $M$
consists of a pair $(g, J)$ where $g$ is a Riemannian
metric and $J:TM\rightarrow TM$ is a complex structure
such that, for any vector fields $X, Y$,
\[\textbf{}
g(X, Y)=g(JX, JY)
\]
and the two--form $\Sigma$ defined by
\[
\Sigma(X, Y)=g(X, JY)
\]
is closed.  
Consider the following local problem
\begin{itemize}
\item Given a Riemannian metric $g$, is there
a K\"ahler metric in the conformal class 
\[[g]=\{cg|c:M\rightarrow \R^+\}\]of $g$?
\end{itemize}
We expect the answer to depend only on the conformal class $[g]$,
thus the local obstructions to the existence of a K\"ahler metric in $[g]$ 
should be given by conformally invariant tensors.

The general Riemannian metric in $2n$ dimensions
depends on $n(2n+1)$ arbitrary functions of $2n$ variables. The 
diffeomorphism freedom can be used to reduce this number
to $2n^2-n$ and finally the freedom of conformally rescaling the metric leaves one with $2n^2-n-1$
arbitrary functions. 
On the other hand the general K\"ahler metric
can be locally described by the K\"ahler potential: there exists a function ${\mathcal K}:M\rightarrow \R$ and
a holomorphic coordinate system $(z^1, \dots, z^n)$ such
that
\[
g=\frac{\p^2{\mathcal K}}{\p z^j\p \ov{z}^k} dz^j d\ov{z}^k.
\]
There is some freedom in the 
K\"ahler potential ${\mathcal K}$, but this freedom 
depends on functions
of $n$ variables. Thus essentially the K\"ahler metric
depends on one arbitrary function of $2n$ variables.
The difference between the number of arbitrary functions
in the general conformal class and the general K\"ahler
metric is
\[
2n^2-n-2
\]
which is positive if $n>1$. This 
gives a lower bound for the number
of conditions a Riemannian metric needs to satisfy in
order to be conformal to a K\"ahler metric. 

If $n=1$ then the general Riemannian metric is conformally flat. In this case one does not need conformal rescaling - any Riemannian metric on a surface is K\"ahler as
in two dimension a conformal structure is equivalent to
a complex structure. 
The first non--trivial case is $n=2$, 
where $(M, g)$ is a real four--manifold. 
There the conformal to K\"ahler condition is 
overdetermined: the numerology tells us to expect
at least four conditions on $g$.

In Section \ref{section1} we derive some necessary conditions on a Riemannian metric in four dimensions 
for it to be locally conformal to K\"ahler. If the self--dual part of the 
conformal curvature of $g$ is non--zero then necessary 
and sufficient conditions are relatively easy to find: The self--dual (SD) Weyl spinor must be of algebraic type $D$ (i.e. it must have two repeated roots when viewed as a symmetric homogeneous polynomial in two variables of degree four), 
and a differential 
obstruction (given in Theorem \ref{non_ASD})  must vanish.

In the anti--self--dual (ASD) case, when the SD Weyl
tensor vanishes, the analysis is more complicated. We shall
construct a natural connection ${\mathcal D}$ on a rank ten vector bundle
\[
E=\Lambda^2_+(M)\oplus \Lambda^1(M)\oplus \Lambda^2_-(M)
\]
over $M$  and 
show (Theorem \ref{theo_tractor}) that there is a one-to-one correspondence between K\"ahler
metrics in the conformal class and non--zero 
parallel sections of ${\mathcal D}$. Readers familiar with the tractor approach to conformal geometry \cite{BEG} 
may want to note that
the vector bundle $E$ 
can be thought of as the space of self--dual
three-forms for the standard tractors $T$ 
of $(M, [g])$, but the induced tractor connection on $(\Lambda^3 T)_+$
is different from ${\mathcal D}$. These two connections
coincide only in the conformally-flat case.

In Section \ref{section_2} we provide examples of metrics which are not conformal to K\"ahler. In the non-ASD case
any metric not of type $D$ gives such an example. In the ASD case we find that an Einstein metric with non--zero cosmological constant can be conformal to 
K\"ahler
if and only if it admits a Killing vector. We then 
argue, using the construction of LeBrun \cite{LeB}, that
ASD Einstein metrics with no symmetries exist. 
Another class of examples is provided by a family
of hyper--hermitian metrics (these are automatically
ASD) which depends on a (locally) harmonic function on the three--sphere 
\cite{Pap,ChValT} . 
In Section \ref{subsection_22} we shall  show that the harmonic function can be chosen so that no K\"ahler metric
exists in the conformal class.

In Section \ref{section_3} we shall make a link between
the `conformal to K\"ahler condition' in dimension four
and the metrisability of projective structures in dimension two. A projective structure on a manifold $U$ is an equivalence class of torsion--free connections which share
the same unparametrised geodesics. In \cite{BDE08}
necessary and sufficient conditions have been determined
in the case when $\dim{U}=2$
for the existence of a (pseudo) Riemannian metric
on $U$ whose geodesics coincide with the geodesics of the given projective structure. If such a metric exists, the projective structure is called metrisable.
Any $n$--dimensional projective structure on $U$ gives rise to a natural conformal structure of signature
$(n, n)$ on $TU$ \cite{walker,YI73}. If $n=2$ this conformal structure
is necessarily ASD \cite{SN,DW07} and we shall show 
(Theorem \ref{theo_projective}) that
the projective structure on a surface $U$ is metrisable if and only if the induced conformal structure on $TU$ admits
a K\"ahler metric or a para--K\"ahler metric. 
This establishes a conjecture made in
\cite{BDE08}. In Section \ref{subsection_31} we shall characterise
the resulting $(2, 2)$ ASD conformal structures
as those which admit a reall parallel section
of the unprimed spin bundle.
\vskip10pt
All considerations in this paper are local. We shall end
this introduction listing some global obstructions\footnote{We thank Claude LeBrun for bringing
these obstructions to our attention.}.
If $(M,[g])$  is a compact ASD conformal Riemannian four-manifold then the existence of a  K\"ahler metric in the
conformal class $[g]$
is equivalent to the following pair of  conditions:
\begin{itemize}
\item  There is a metric in $[g]$ with vanishing scalar curvature.
\item The number $b_+(M)$ of positive eigenvalues of the intersection
form  is non-zero.
\end{itemize}
The second condition is purely topological, 
while the first condition is equivalent to vanishing of the
Yamabe constant of $[g]$. These conditions impose
constraints on the diffeomorphism type of $M$. 
The only simply connected $M$ that are allowed are
$K3$ and $\CP_2 \# k\overline{\CP}_2$, $k\geq 10$.
See \cite{Lebrun86, KLP97,  RS05, LM08} for details.

\vspace{2ex}{\bf Acknowledgements.} 
We are grateful to  Mike Eastwood and 
Claude LeBrun for helpful discussions. 
We also thank the anonymous referee for suggesting 
the example given in Section \ref{referee_example} to us.

MD was partly supported by the MISGAM 
and ENIGMA programs of the European Science Foundation.
\subsection{Spinors in four dimensions}
Given an oriented Riemannian four-manifold $(M,g)$, the
Hodge operator $\ast:\Lambda^2\rightarrow \Lambda^2$
satisfies $\ast^2=\mbox{Id}$
and induces a decomposition
\be 
\label{splitting}
\Lambda^{2} = \Lambda_{+}^{2} \oplus \Lambda_{-}^{2}
\ee
of 2-forms
into self--dual
and anti--self--dual\index{anti--self--dual} components. 
The rank three vector bundles $\Lambda_\pm^2\rightarrow M$ 
are the eigenspaces of $\ast$ with eigenvalues $\pm 1$ respectively.
The Riemann tensor
has the index symmetry $R_{abcd}=R_{[ab][cd]}$ so can be thought of
as a map $\mathcal{R}: \Lambda^{2} \rightarrow \Lambda^{2}$.
This map decomposes under (\ref{splitting}) as follows:
\be \label{decomp}
{\mathcal R}=
\left(
\mbox{
\begin{tabular}{c|c}
&\\
$C_++\frac{R}{12}$&$\phi$\\ &\\
\cline{1-2}&\\
$\phi$ & $C_-+\frac{R}{12}$\\&\\
\end{tabular}
} \right) .
\ee
The $C_{\pm}$ terms are the self-dual and anti-self-dual parts
of the Weyl tensor, the $\phi$ terms are the
trace-free Ricci curvature, and $R$ is the scalar curvature which acts
by scalar multiplication. The Weyl tensor\index{Weyl tensor} is  conformally
invariant, so can be
thought of as being defined by the conformal structure 
$[g]$.

Locally there exist complex rank two vector bundles
$\spp, \spp'$  (called spin-bundles) over $M$
equipped with parallel symplectic
structures
$\ep, \ep'$ such that
\be
\label{can_bun_iso}
\C\otimes T M\cong {\spp}\otimes {\spp'}
\ee
is a  canonical bundle isomorphism, and
\[
g(v_1\otimes w_1,v_2\otimes w_2)
=\varepsilon(v_1,v_2)\varepsilon'(w_1, w_2)
\]
for $v_1, v_2\in \Gamma(\spp)$ and $w_1, w_2\in \Gamma(\spp')$.
We use the standard convention \cite{PR} in which spinor indices are capital
letters, unprimed for sections of $\spp$ and primed for sections of $\spp'$. For example $\mu_{A}$ denotes a section of $\spp^{*}$, the dual of $\spp$, and $\nu^{A'}$ a section of $\spp'$.
The symplectic structures $\varepsilon_{AB}$
and $\varepsilon_{A'B'}$ (such that
$\varepsilon_{01}=\varepsilon_{0'1'}=1$) are used to raise and
lower indices by $\mu_{A}:=\mu^{B}\varepsilon_{BA}, 
\mu^A=\varepsilon^{AB}\mu_B$. 

The decomposition (\ref{splitting})
of  two--forms takes a simple form 
in the spinor notation. If 
$
F_{ab}=F_{[ab]}
$ is a
two--form
then
\be 
\label{2form}
F_{AA'BB'} = f_{AB}\ep_{A'B'} + \tilde{f}_{A'B'}\ep_{AB},
\ee
where $f_{AB}$ and $\tilde{f}_{A'B'}$ are symmetric in their indices.  This is precisely
the decomposition of $F$ into self-dual and anti-self dual parts.
Thus
\[
\Lambda^{2}_{+} \cong \spp'^{*} \odot \spp'^{*}, \qquad
\Lambda^{2}_{-} \cong \spp^{*} \odot\spp ^{*}.
\]

In terms of the spinor notation the decomposition 
(\ref{decomp})
of the Riemann tensor is
\begin{eqnarray} \label{riemann}
R_{abcd} &=& \psi_{ABCD} \ep_{A'B'} \ep_{C'D'} + {\psi}_{A'B'C'D'}
\ep_{AB} \ep_{CD} \nonumber\\ 
&&+ \phi_{ABC'D'} \ep_{A'B'}\ep_{CD} +
\phi_{A'B'CD} \ep_{AB}\ep_{C'D'} \nonumber\\&& + \frac{R}{12} (\ep_{AC}\ep_{BD}
\ep_{A'C'} \ep_{B'D'} - \ep_{AD} \ep_{BC} \ep_{A'D'} \ep_{B'C'}),
\end{eqnarray}
where $\psi_{ABCD}$ and $\psi_{A'B'C'D'}$ are ASD and SD Weyl spinors
which are symmetric in their indices and 
$\phi_{A'B'CD}=\phi_{(A'B')(CD)}$ is the traceless Ricci spinor.
A conformal structure is called ASD iff $\psi_{A'B'C'D'}=0$.
\section{Conformal-to-K\"ahler conditions}
\label{section1}
Let $(M, g)$ be a four manifold such that 
\be
\label{conf_resc_form}
\hat{g}_{ab}=\Omega^2 g_{ab}
\ee
is a K\"ahler metric
with a K\"ahler form $\hat{\Sigma}$.
The two--form $\hat{\Sigma}$ on $M$ induces a natural orientation given by the volume form $\hat{\Sigma}\wedge\hat{\Sigma}$.
With respect to this orientation 
$\hat{\Sigma}$ is self--dual. Thus, 
in terms of 2-component spinors, the K\"ahler form may be written as
\[\hS_{ab}=\hw_{A'B'}\hat{\varepsilon}_{AB},\]
so that
\begin{eqnarray*}
\hnabla_{AA'}\hw_{B'C'}&=&0,\\
\hw_{A'B'}\hw^{A'B'}&=&2.
\end{eqnarray*}
Set
\[
\Upsilon_{a}=\Omega^{-1}\nabla_a\Omega.
\]
Under conformal rescaling we shall make the choice that $\hw_{A'B'}$ transforms as
\be
\hw_{A'B'}=\Omega^2\omega_{A'B'},
\ee
since then
\be\label{res1}
\hnabla_{AA'}\hw_{B'C'}=\Omega^{2}\left(\nabla_{AA'}\omega_{B'C'}-
\Upsilon_{B'A}\omega_{A'C'}-\Upsilon_{C'A}\omega_{B'A'}
+2\Upsilon_{AA'}\omega_{B'C'}\right)=0,
\ee
so that
\[
\hnabla_{A(A'}\hw_{B'C')}=\Omega^{2}\nabla_{A(A'}\omega_{B'C')}.\]
We have therefore arrived at the following result of Pontecorvo 
\cite{pontecorvo}.
\begin{lemma}
\label{first_lemma}
The metric $g_{ab}$ is conformal to a K\"ahler 
metric if and only if there exists a real, symmetric spinor field $\omega_{A'B'}$ satisfying
\be\label{ck0}
\nabla_{A(A'}\omega_{B'C')}=0,
\ee
and such that $\omega_{A'B'}\omega^{A'B'}\neq 0$.
\end{lemma}
{\bf Proof.}
The condition (\ref{ck0}) is equivalent to 
\be\label{ck1}
\nabla_{AA'}\omega_{B'C'}=\ep_{A'B'}K_{C'A}+\ep_{A'C'}K_{B'A}.
\ee
for some one--form $K_a$.
Note that (\ref{ck1}) follows from (\ref{res1}) with
\be\label{om0}
K_{AA'}=-\omega_{A'}^{\;\;B'}\Upsilon_{B'A},
\ee
and $\Omega$ can be found from
\be\label{om1}
\omega_{A'B'}\omega^{A'B'}=2\Omega^{-2}.
\ee
Suppose instead we have a real solution of (\ref{ck0}), so that also (\ref{ck1}) holds for some $K_{AA'}$. Define $\Omega$ by (\ref{om1}), then one may calculate that
\be\label{om3}
\Upsilon_{AA'}=\Omega^2\omega_{A'}^{\;\;C'}K_{C'A},
\ee
and then the right-hand-side of (\ref{res1}) is zero, so that
$g_{ab}$ is indeed conformal to K\"ahler, provided $\Omega$ has no zeroes. Thus the existence of a real solution of (\ref{ck0}) is necessary and sufficient for the metric to be conformal to K\"ahler, in regions where $\Omega\neq 0$.\koniec

We want to find conditions on curvature for a suitable solution of (\ref{ck1}) to exist. We shall use the procedure of prolongation
to introduce new variables, and rewrite (\ref{ck0}) as a closed overdetermined system of first order PDEs
which then can be investigated for Frobenius integrability.

Differentiate and commute derivatives on (\ref{ck1}) to obtain
\be\label{notSD}
\psi^{E'}_{\;\;(A'B'C'}\omega_{D')E'}=0, \ee 
in terms of the SD Weyl spinor $\psi_{E'A'B'C'}$, and 
\be\label{ck2}
\nabla_{AA'}K_{BB'}=-P_{ABA'C'}\omega_{B'}^{\;\;C'}+\ep_{A'B'}\rho_{AB}
\ee
in terms of some as yet unknown but symmetric $\rho_{AB}$, where $P_{ABA'B'}$ is given in terms of the Ricci spinor $\phi_{ABA'B'}$ and curvature scalar $\Lambda=R/24$ by
\[
P_{ABA'B'}=\phi_{ABA'B'}-\Lambda\ep_{AB}\ep_{A'B'}.
\]

\subsection{The non anti--self--dual case}
Suppose that $\psi_{A'B'C'D'}\neq 0$, then (\ref{notSD})
has no solution for non-null (hence for real) $\omega_{A'B'}$ unless
\be
\label{notSD2}
\psi_{A'B'C'D'}=c\Theta_{(A'B'}\Theta_{C'D')},
\ee 
for
some function $c$ and real symmetric $\Theta_{A'B'}$
(one can always choose $\Theta_{A'B'}$ to set $c$ to $\pm 1$). 
If
(\ref{notSD2}) does hold, which is equivalent to type-D SD Weyl spinor in the Petrov-Pirani-Penrose classification (
see e.g. \cite{PR}, or \cite{Der83,AG97} for discussion relevant in 
the compact case), then any
solution $\omega_{A'B'}$ of (\ref{notSD}) is proportional to
$\Theta_{A'B'}$. Thus in the case $\psi_{A'B'C'D'}\neq
0$, one may determine whether the metric is conformal to K\"ahler by simply checking if any candidate $\omega_{A'B'}=f\Theta_{A'B'}$ solves (\ref{ck1})  and is appropriately nonzero.
Following \cite{PR} define two conformal invariants 
\[
{\mathcal I} = |\psi|^2=\psi_{A'B'C'D'} \psi^{A'B'C'D'}, \quad {\mathcal J} 
= \psi_{A'B'}^{\;\;\;\;\;\;\;\; C'D'}\psi_{C'D'}^{\;\;\;\;\;\;\;\; E'F'}
\psi_{E'F'}^{\;\;\;\;\;\;\;\; A'B'}
\]
of conformal weight $-4$ and $-6$ respectively.
Define a spinor
\be
\label{defi_of_t}
T_{AA'B'C'D'E'}:=\nabla_{A(A'}\psi_{B'C'D'E')}-
V_{A(A'}\psi_{B'C'D'E')}, 
\ee
where
\be
\label{exp_for_V}
V_{AA'}=\frac{1}{|\psi|^2}\Big(\frac{1}{6}\nabla_{AA'}|\psi|^2+
\frac{4}{3}\psi^{B'C'D'E'}\nabla_{AB'}\psi_{C'D'E'A'}\Big).
\ee
\begin{theo}
\label{non_ASD}
Let $(M, [g])$ be a conformal manifold such that the self--dual
part of the conformal curvature of $[g]$ is non--zero. Then
there exists a K\"ahler metric in the conformal class if
and only if the following conformally invariant
conditions hold
\be
\label{alg_cond_IJ}
{\mathcal I}^3=6{\mathcal J}^2
\ee
\be
\label{T_condition}
T_{AA'B'C'D'E'}=0,
\ee
\be
\label{closure_of_V}
\nabla_{[a} V_{b]}=0.
\ee
\end{theo}
{\bf Proof.}
Consider a metric $g$ with a SD Weyl tensor satisfying
(\ref{notSD2}). The Weyl spinor corresponds to a homogeneous
polynomial in two variables of degree 4 with two repeated roots,
and (\ref{alg_cond_IJ}) holds identically. If the conformal class
of $g$ contains a K\"ahler metric then $\Theta_{A'B'}$
satisfies (\ref{ck0}). Differentiating (\ref{notSD2}) and
symmetrising over the primed indices gives
\be
\label{non_SD1}
\nabla_{A(A'}\psi_{B'C'D'E')}=V_{A(A'}\psi_{B'C'D'E')}
\ee
where $V_{AA'}$ is a gradient. We shall analyse (\ref{non_SD1})
for general $V_{AA'}$. Contracting both sides
with $\psi^{B'C'D'E'}$ and using the identity
\[
\psi^{B'C'D'E'}\psi_{B'C'D'A'}=\frac{1}{2}\delta^{E'}_{A'}|\psi|^2
\]
yields the expression (\ref{exp_for_V}) for $V$.
Substituting this back to (\ref{non_SD1}) gives the vanishing
of $T_{AA'B'C'D'}$ defined by (\ref{defi_of_t}). Now we need
to establish the conformal invariance of conditions
(\ref{T_condition}) and  (\ref{closure_of_V}). 
Under the conformal rescaling
(\ref{conf_resc_form}) we have
\[
\hat{\varepsilon}_{A'B'}= \Omega\varepsilon_{A'B'}, \quad
\hat{\varepsilon}^{A'B'}= \Omega^{-1}\varepsilon^{A'B'}
\]
and
\[
\hat{\psi}_{A'B'C'D'}= \psi_{A'B'C'D'}, \quad
\hat{\psi}^{A'B'C'D'}= \Omega^{-4}\psi^{A'B'C'D'}.
\]
Using
\begin{eqnarray*}
\hat{\nabla}_{AB'}\hat{\psi}_{C'D'E'A'}&=&
\nabla_{AB'}{\psi}_{C'D'E'A'}\\
&&-\Upsilon_{AC'}{\psi}_{B'D'E'A'}
-\Upsilon_{AD'}{\psi}_{B'C'E'A'}
-\Upsilon_{AE'}{\psi}_{B'C'D'A'}
-\Upsilon_{AA'}{\psi}_{B'C'D'E'}
\end{eqnarray*}
gives
\be
\label{confo_V}
\hat{V}_{AA'}=V_{AA'}-4\Upsilon_{AA'}
\ee
and
\[
\hat{\nabla}_{A(B'}\hat{\psi}_{C'D'E'A')}=
{\nabla}_{A(B'}{\psi}_{C'D'E'A')}-4V_{A(B'}{\psi}_{C'D'E'A')}.
\]
Therefore 
\[
\hat{T}_{AA'B'C'D'}={T}_{AA'B'C'D'}
\]
and condition (\ref{T_condition}) is conformally invariant.
The closure of $V$ is also conformally invariant because
of (\ref{confo_V}).

Conversely, assume that conditions
(\ref{alg_cond_IJ}), (\ref{T_condition})
and 
(\ref{closure_of_V}) hold. The algebraic condition 
$(\ref{alg_cond_IJ})$ implies that $\psi_{A'B'C'D'}$ has a repeated root.
There are various possibilities, but the only one compatible
with the Riemannian signature of $g$ is that $\psi_{A'B'C'D'}$ 
is given by
(\ref{notSD2}) for some $\Theta_{A'B'}$. Now construct
$V_{AA'}$ and $T_{AA'B'C'D'}$. The closure condition
(\ref{closure_of_V}) implies that locally $V$ is a gradient.
The relation
(\ref{confo_V}) can now be used to find a conformal factor
such that $V=0$ and (\ref{T_condition}) reduces
to
\[
{\nabla}_{A(A'}{\psi}_{B'C'D'E')}=0
\]
which implies
\[
\Theta_{(B'C'}{\nabla^A}_{A'}\Theta_{D'E')}=0
\]
where we have reabsorbed $\phi$ into $\Theta$. This is equivalent to
\[
\nabla_{A(A'}\Theta_{B'C')}=0.
\]
Thus, by Lemma \ref{first_lemma}, there exists a K\"ahler metric in the conformal class of $g_{ab}$.
\koniec

\subsection{The anti--self--dual case}
The harder case is therefore zero SD Weyl spinor, so henceforth we shall assume
that $\psi_{A'B'C'D'} = 0$ but $\psi_{ABCD}\neq 0$ i.e. the Weyl tensor is nonzero but ASD.

Commute derivatives on (\ref{ck2}) to obtain
\be\label{ck3}
\nabla_{AA'}\rho_{BC}=\omega_{A'}^{\;\;E'}\nabla_{E'}^{\;\;D}\psi_{ABCD}-K_{A'}^{\;\;D}\psi_{ABCD}-2P_{A'E'A(B}K_{C)}^{\;\;E'}
\ee
using curvature assumptions already made.

Suppose for a moment that $g_{ab}$ actually is K\"ahler. Then (\ref{ck0}, \ref{ck1}, \ref{ck2}) must hold with constant $\Omega$, in which case $K_a$ vanishes but from (\ref{ck2}) we deduce
\[\Lambda=0\;;\;\phi_{ABA'B'}=\omega_{A'B'}\rho_{AB}\]
and we may identify $\rho_{AB}\epsilon_{A'B'}$ as the Ricci form. In this case (\ref{ck3}) becomes an identity and this gives an understanding of what $\rho$ is in general. 

In the general case, we can next commute derivatives on $\rho$. In terms of $\Delta_{A'B'}:=\nabla_{A(A'}\nabla_{B')}^{\;\;A}$, we obtain an identity from $\Delta_{A'B'}\rho_{AB}$ (using the vanishing of the Bach tensor, which holds for any metric with ASD or SD Weyl tensor \cite{PR}) but something new arises from $\Delta_{AB}\rho_{CD}$ namely the necessary condition:

\bea\label{cond}
&&4\psi_{E(ABC}\rho_{D)}^{\;\;E}+K^{EE'}\nabla_{EE'}\psi_{ABCD}+
4K_{E'(A}\nabla^{EE'}\psi_{BCD)E}\\&&+\omega_{E'A'}\left(\nabla_{(A}^{\;\;A'}\nabla^{EE'}\psi_{BCD)E}+\psi_{E(ABC}\phi_{D)}^{\;\;EE'A'}\right)=0.\nonumber
\eea
This represents five linear restrictions on the ten-component column vector
\[X_{\mathcal{A}}=(\omega_{A'B'},K_{AA'}, \rho_{AB})^T.\] 
This is not enough to be overdetermined, but it and its derivatives have to hold, which eventually gives constraints on the conformal curvature, as we shall see with an example below.

\medskip

\subsection{The prolongation bundle}
\label{subsection_ASD}
The formalism has led us to consider a natural rank-10 vector bundle\footnote{A similar construction, albeit not in a conformally invariant setting, was used in \cite{S03} 
in the study of  conformal Killing forms.}
\[
E=\Lambda^2_+(M)\oplus \Lambda^1(M)\oplus \Lambda^2_-(M)
\]
with sections $X_{\mathcal{A}}$ which have the following behaviour under conformal rescaling of the metric:
\be
\label{10-tractor}
\hat{X}_{\mathcal{A}}:=
\left(\begin{array}{c}\hw_{A'B'}\\\hk_{AA'}\\\hrho_{AB}
\end{array} \right) 
=
\left(\begin{array}{c}\Omega^2\omega_{A'B'}\\\Omega(K_{AA'}+\Upsilon_{AB'}\omega_{A'}^{\;\;B'})\\\rho_{AB}+2\Upsilon_{C'(A}K_{B)}^{\;\;C'}+\omega^{C'D'}\Upsilon_{C'A}\Upsilon_{D'B}
\end{array} \right) 
\ee

This transformation recalls that of a tractor, \cite{BEG}, and as we noted in the introduction the bundle $E$ can be understood as the bundle of self-dual 3-forms $\Lambda^3_+(T)$ where $T$ is the usual 
tractor bundle (although with a different connection). Define a derivative $\mathcal{D}$ on this vector bundle by
\be
\label{tractorD}
{\mathcal{D}}_a{X}_{\mathcal{B}}:=
\left(\begin{array}{c}
\nabla_{AA'}\omega_{B'C'}-\ep_{A'B'}K_{C'A}
-\ep_{A'C'}K_{B'A}\\
\nabla_{AA'}K_{BB'}+P_{ABA'C'}\omega_{B'}^{\;\;C'}-\ep_{A'B'}\rho_{AB}\\
\nabla_{AA'}\rho_{BC}-\omega_{A'}^{\;\;E'}\nabla_{E'}^{\;\;D}\psi_{ABCD}+K_{A'}^{\;\;D}\psi_{ABCD}+2P_{A'E'A(B}K_{C)}^{\;\;E'}
\end{array} \right). 
\ee
We can now state our result in the following way
\begin{theo}
\label{theo_tractor}
A four--dimensional Riemannian manifold $(M, g)$
with ASD conformal curvature 
is locally conformal to K\"ahler if and only if there exists a parallel non--zero section ${X}_{\mathcal{B}}$ for $\mathcal{D}$:
\be\label{cond3}{\mathcal{D}}_a{X}_{\mathcal{B}}=0\ee
such that $\omega_{A'B'}\omega^{A'B'}\neq 0$.
\end{theo}
We calculate the curvature of this connection from the general formula
\[{\mathcal{D}}_{[a}{\mathcal{D}}_{b]}{X}_{\mathcal{C}}=\frac12{\mathcal{R}}_{ab\,{\mathcal{C}}}^{\;\;\;\;\;{\mathcal{E}}}X_{\mathcal{E}},\]
to find that ${\mathcal{R}}_{ab\,{\mathcal{C}}}^{\;\;\;\;\;\;\;{\mathcal{E}}}$ is ASD
on the index pair $ab$, so that the only nonzero components are ${\mathcal{R}}_{AB{\mathcal{C}}}^{\;\;\;\;\;\;\;{\mathcal{E}}}$.

From (\ref{cond3}) we may deduce a necessary condition to be conformal to K\"ahler as
\[{\mathcal{R}}_{ab\,{\mathcal{C}}}^{\;\;\;\;\;{\mathcal{E}}}X_{\mathcal{E}}=0,\]
 and we claim that
(\ref{cond}) is precisely
\[{\mathcal{R}}_{AB\,{\mathcal{C}}}^{\;\;\;\;\;\;{\mathcal{E}}}X_{\mathcal{E}}=0.\]

To see how many conditions follow by differentiating (\ref{cond}), we calculate the
Bianchi identity:
\[{\mathcal{D}}_{[a}{\mathcal{R}}_{bc]\,{\mathcal{C}}}^{\;\;\;\;\;\;\;{\mathcal{E}}}=0\]
which given ASD reduces to
\[{\mathcal{D}}_{A'}^{\;A}{\mathcal{R}}_{AB \,{\mathcal{C}}}^{\;\;\;\;\;\;\;\;{\mathcal{E}}}=0.\]
Thus from differentiating (\ref{cond}), we obtain just 
\be\label{cond4}
{\mathcal{D}}_{A'(A}{\mathcal{R}}_{BC)\,{\mathcal{C}}}^{\;\;\;\;\;\;\;\;{\mathcal{E}}}X_{\mathcal{E}} =0.
\ee
This turns out to be
a further twelve conditions, so that we have seventeen 
linear conditions on the ten-component $X_{\mathcal{A}}$ .

\medskip

\section{Examples of Riemannian four-metrics not conformal to K\"ahler}
\label{section_2}

\subsection{Nonzero SD Weyl spinor}

A Riemannian 4-metric whose Weyl tensor has a non-zero SD Weyl spinor which is not type D cannot be conformal to K\"ahler. Examples could be found among the self--dual vacuum metrics, 
though these would be K\"ahler with the opposite orientation. 

Other examples may be given as follows: consider the metric
\[g=(\delta_{ab}+\frac12R^0_{acbd}x^cx^d)dx^adx^b,\]
where $R^0_{acbd}$ is a constant tensor with Riemann tensor symmetries. The Riemann tensor of this metric, calculated at the origin of the coordinates $x^a$, is $R_{acbd}(0)=R^0_{acbd}$ so that, if this does not have type-D SD Weyl tensor then this metric is not conformal to K\"ahler in any neighbourhood of the origin. By including more terms in the Taylor series, one may construct examples which satisfy the first condition of Theorem \ref{non_ASD} but not the 
later ones.

\subsection{ASD Einstein} 

Examples of ASD conformal classes which do not contain a K\"ahler
metric would be given by ASD Einstein metrics with non--zero
cosmological constant, which did not admit any Killing vectors.
This is because of
\begin{prop}
An anti-self--dual 
Einstein metric $g$ with $\Lambda\neq 0$ is conformal to a K\"ahler
metric iff $g$ admits a Killing vector. 
\end{prop}
{\bf Proof.}
Let $g$ be ASD Einstein, and  assume that there exists
a K\"ahler metric in $[g]$.
Then (\ref{ck1}) and (\ref{ck2}) imply that
$K=0$ or $K$ is Killing. 

If $K=0$ then $\omega$ is parallel so
$g$ is already K\"ahler--Einstein--ASD. But K\"ahler--ASD implies
scalar--flat. Thus $g$ is Ricci--flat  K\"ahler and
so hyper--K\"ahler.

If $K$ is Killing then one can locally 
choose coordinates such that $K=\p/\p t$ and
$g$ is given by \cite{Prz}, \cite{T97}
\[
g=\frac{P}{z^2}(e^u(dx^2+dy^2)+dz^2)+\frac{1}{Pz^2}(d t +\alpha)^2
\]
where $u=u(x, y, z)$  is a solution of the $SU(\infty)$ Toda equation
\[
u_{xx}+u_{yy}+({e^u})_{zz}=0,
\]
the function $P$ is given by $2\Lambda P=zu_z-2$, and
\[
d\alpha=-P_xdy\wedge dz-P_y dz\wedge dx -(Pe^u)_zdx\wedge dy.
\]
Thus $\hat{g}=z^2g$ is of the form given by LeBrun ansatz \cite{L91}
because $P$ satisfies the linearised  
$SU(\infty)$ Toda equation. Therefore 
$\hat{g}$ is K\"ahler with vanishing Ricci scalar. 
All such metrics have ASD conformal curvature.
\koniec
It is worth pointing out that such an ASD Einstein metric may be conformal to K\"ahler in more than one way. Thus for example 
$\CP^2$ with the Fubini-Study metric, but the wrong orientation so that it is not K\"ahler, is ASD Einstein and is conformal to K\"ahler in 
different ways via the above construction using the Killing 
vectors $\p/\p\psi$ and $\p/\p\phi$ in the usual Euler angles.

To find ASD Einstein metrics with negative cosmological constant not conformal to K\"ahler, one could apply the LeBrun construction,  \cite{LeB}, which fills in a collar neighbourhood of a compact metric at infinity with such a metric. Now if one takes a metric on $S^3$ with no symmetries as the metric at infinity, then the resulting Einstein metric will also have no symmetries and so cannot be conformal to K\"ahler. 

\subsection{Compact hyperbolic four--manifolds}
\label{referee_example}

This construction, suggested to us by the referee, gives  examples of  
Riemannian manifolds which are locally conformal to K\"ahler in many ways,
but are not globally conformally K\"ahler. Consider a compact hyperbolic 
four manifold $(M, g)$. Its traceless Ricci and Weyl tensors vanish 
and the Ricci scalar is equal to $-1$. In  particular the Weyl curvature of $g$ is ASD, so
$g$ is not K\"ahler as its Ricci scalar is not zero. Let us assume that there 
exists globally defined non-degenerate $\omega_{A'B'}$ 
such that (\ref{ck0}) holds.
The spinor $\omega_{A'B'}$ is not parallel (otherwise $g$ would be K\"ahler),
thus (\ref{ck1}) defines a non--zero one--form $K_{AA'}$ 
and formula (\ref{ck2})
implies that $K^{AA'}$ is a Killing vector. However this is not possible as
the Bochner argument implies that non--trivial Killing vector fields
can not exist on compact manifold with negative Ricci curvature: A 
Killing vector field $K^a$ satisfies $\Box K_a+R_{ab}K^b=0$, where in
our case $R_{ab}=-g_{ab}/4$. Multiplying this identity by $K^a$ and
integrating over a compact manifold $M$ yields
\[
\int_M |\nabla K|^2 \sqrt{g}d^4 x =-\frac{1}{4}\int_M |K|^2 \sqrt{g}d^4 x.
\]
Thus $K^a$ must vanish everywhere which contradicts our assumption about
$\omega_{A'B'}$.

\subsection{Hypercomplex examples}
\label{subsection_22}
Families of these were given by \cite{Pap} and \cite{ChValT}. The metric takes the form:

\be
\label{met3}
g=V^2(\sigma^{\;2}_1+\sigma^{\;2}_2+\sigma^{\;2}_3)+V^{-1}(dt+\alpha)^2,
\ee
where the basis $(\sigma_i)$ of invariant one-forms on $S^3$ satisfies
\[d\sigma_1=\sigma_2\wedge\sigma_3\;\mathrm{   etc}\]
and $V$ and $\a$ satisfy
\[d\alpha=V_1\sigma_2\wedge\sigma_3+V_2\sigma_3\wedge\sigma_1
+V_3\sigma_1\wedge\sigma_2\]
where we define $V_i$
 by
\[dV=V_1\sigma_1+V_2\sigma_2+V_3\sigma_3.\]
Necessarily then 
\[\Delta V:=V_{11}+V_{22}+V_{33}=0,\]
so that $V$ is harmonic on $S^3$. For a metric which isn't conformally flat, $V$ isn't constant and so necessarily has singularities, so that these metrics are only defined locally. It is straightforward to check that the Weyl tensor is ASD and the Killing vector $T=\partial/\partial t$ has ASD derivative. 

For the examples, the idea is to choose $V$ so that, at a point $p\in S^3$, $V=1$ while the 
first and second derivatives $V_i$ and $V_{ij}$ all vanish. By considering harmonic functions expanded in spherical harmonics on $S^3$ it is clear that this can be done and then 
the third partials $V_{ijk}(p)$ are freely disposable, subject to being symmetric and 
trace-free; that gives a 7-dimensional vector space of choices for $V_{ijk}(p)$ (but an infinite-dimensional family of choices for $V$).

For the derivative of the Killing vector, we calculate
\[\nabla_aT_b=\varepsilon_{A'B'}\phi_{AB}\]
with
\[\phi_{AB}=\frac1VT_{A}^{\;B'}\nabla_{BB'}V.\]
From the Killing vector identity 
\[\nabla_a\nabla_bT_c=R_{bcad}T^d\]
we may calculate
\[\psi_{ABCD}=2VT_{D}^{\;C'}\nabla_{C'A}\phi_{BC}.\]
Thus, at $p$,
\[ \psi_{ABCD}=0= \nabla_aT_b=\nabla_a\nabla_bT_c.\]
We need explicit formulae for the Ricci curvature, for which we use Cartan calculus to find
\be
\label{f1}
R_{ab}=\frac12(T_aT_b-\frac1Vg_{ab})
\ee
so that, at $p$,
\[ \nabla_aR_{bc}=0=\nabla_a\nabla_bR_{cd}.\]
From the Bianchi identity this means that, at $p$,
\[\nabla^{AA'} \psi_{ABCD}=0= \nabla^{FF'}\nabla^{AA'} \psi_{ABCD}.\]
The Lie derivative of the Weyl curvature along the Killing vector necessarily vanishes everywhere, so that, at $p$,
\be\label{f2}
T^{EE'}\nabla_{EE'} \psi_{ABCD}=0= T^{EE'}\nabla_{FF'}\nabla_{EE'} \psi_{ABCD}.
\ee
If we use all that we have so far established in (\ref{cond}) evaluated at $p$, then it collapses to
\be
\label{k1}
K^{EE'}\nabla_{EE'} \psi_{ABCD}=0.
\ee
Clearly this is satisfied with $K^a(p)$ any multiple of the Killing vector $T^a(p)$; we claim that generically there are no other solutions. To prove this we need explicit formulae for the components of the Weyl tensor. These are readily found via Cartan calculus in the tetrad:
\[\theta^0=V^{-1/2}(dt+\alpha),\;\theta^1=V^{1/2}\sigma_1, \;\theta^2=V^{1/2}\sigma_2,\;\theta^3=V^{1/2}\sigma_3.\] 
We find
\bea\label{weyl2}
C_{0101}&=&-\frac{V_{11}}{2V^2}-\frac{1}{2V^3}(V_2^{\;2}+V_3^{\;2}-2V_1^{\;2})\\
C_{0102}&=&-\frac{V_{12}}{2V^2}-\frac{3}{2V^3}V_1V_2-\frac{V_3}{4V^2}
\eea
from which the rest follow by permutations and anti-self-duality. Now it follows that, at $p$, the components of the derivative $\nabla_aC_{bcde}$ are proportional to the third derivatives $V_{ijk}$. Condition (\ref{k1}) is equivalent to
\[K^a\nabla_aC_{bcde}=0\] 
which reduces to an equation on the components of $K^a$ orthogonal to $T^a$:
\be
\label{co2}
K^iV_{ijk}=0
\ee
and this implies $K^i=0$ for generic $V_{ijk}$ (and in particular for the choice we make below). Thus
\[K^a(p)=c_1T^a(p),\]
for some real constant $c_1$, which could be zero. 
We move on to the derivative of (\ref{cond}), evaluated at $p$. This is messy but straightforward. We are assisted by an identity obtained from the Bianchi identity
\bea
\nabla^{AA'} \psi_{ABCD}&=&\nabla^{B'}_{(B} \phi_{CD)B'}^{\;\;\;\;\;\;\;\;\;\;A'}\\
&=&-\frac{3}{4}\phi_{(AB}T_{C)}^{\;\;A'}
\eea
by using (\ref{f1}) for Ricci. This simplifies the last term in (\ref{cond}) to
\[\omega_{E'A'}T^{E'}_{(A}T^{A'E}\psi_{BCD)E} \]
in which form its derivative is easier to see. We take the derivative $\nabla_F^{F'}$ of (\ref{cond}), and evaluate at $p$, using $K^a(p)=c_1T^a(p)$. The resulting expression is simplified by defining
\[\hat\rho_{AB}=\rho_{AB}-\frac14T^{E'}_{A}T^{F'}_{B}\omega_{E'F'}\]
when it becomes
\[5\hat\rho_{(A}^{E}\nabla_{B}^{F'} \psi_{CDF)E}+\frac18T^{E}_{A'}\omega^{F'A'}T_{G'F}\nabla_{E}^{G'} \psi_{ABCD} -\Lambda\omega^{E'F'}\nabla_{E'F}\psi_{ABCD}=0.\]
At $p$, $V=1$ so that $T_aT^a=1$ and, by (\ref{f1}), $\Lambda=-1/16$. The second term in the above expression simplifies further using (\ref{f2}) and some algebra as follows:
\[\frac18T^{E}_{A'}\omega^{F'A'}T_{G'F}\nabla_{E}^{G'} \psi_{ABCD}=\frac{1}{16}\omega^{F'A'}\nabla_{FA'}\psi_{ABCD},\]
which duplicates the third term. We finally obtain
\be
\label{cond2}
5\hat\rho_{(A}^{E}\nabla_{B}^{F'} \psi_{CDF)E}+\frac{1}{8}\omega^{E'F'}\nabla_{E'F}\psi_{ABCD}=0.
\ee

Note that (\ref{cond2}) is totally symmetric on the indices $ABCDF$ (as expected - see (\ref{cond4})) so it is a system of twelve (real) linear equations on the unknowns $(\omega_{A'B'},\rho_{AB})$. If the relevant determinant is non-zero then the only solution will be
\[\omega_{A'B'}=0=\rho_{AB}
\]
holding at $p$. By (\ref{om1}), a zero in $\omega_{A'B'}$ corresponds to a pole in $\Omega$, which of course isn't allowed. This will show that this metric is not conformally K\"ahler in any neighbourhood of $p$. Thus we must analyse (\ref{cond2}) further.

We introduce two spinor dyads $(o^A,o^{\dagger A})$ and $(\tilde{o}^{A'},\tilde{o}^{\dagger A'})$ normalised by
\[o_Ao^{\dagger A}=1=\tilde{o}_{A'}\tilde{o}^{\dagger A'}\] and related to the tetrad by
\begin{eqnarray*}
\frac{1}{\sqrt{2}}(\theta^0+i\theta^1)&=&o_A\tilde{o}_{A'}dx^a\\
\frac{1}{\sqrt{2}}(\theta^0-i\theta^1)&=&o^{\dagger}_A\tilde{o}^\dagger_{A'}dx^a\\
\frac{1}{\sqrt{2}}(\theta^2+i\theta^3)&=&o^{\dagger}_A\tilde{o}_{A'}dx^a\\
\frac{1}{\sqrt{2}}(\theta^2-i\theta^3)&=&-o_A\tilde{o}^{\dagger}_{A'}dx^a,
\end{eqnarray*}
with the corresponding operators:
\begin{eqnarray*}
D&=&\frac{1}{\sqrt{2}}(e_0+ie_1)\\
\overline{D}&=&\frac{1}{\sqrt{2}}(e_0-ie_1)\\
\delta&=&\frac{1}{\sqrt{2}}(e_2+ie_3)\\
\overline\delta&=&\frac{1}{\sqrt{2}}(e_2-ie_3).
\end{eqnarray*}
We introduce a spinor field $\chi_{ABCDE}$ by
\[\nabla_{E'E}\psi_{ABCD}=\tilde{o}_{E'}\chi_{ABCDE}+\;\;\mbox{H.c.}
\]
where, at $p$, $\chi_{ABCDE}$ is totally symmetric (though not Hermitian). For (\ref{cond2}) we need the dyad components of $\chi_{ABCDE}$. In the 
Newman--Penrose  formalism \cite{PR} at $p$ these can be written
\begin{eqnarray*}
\chi_0&=&\overline\delta^3V\\
\chi_1&=& -\overline{D}\, \overline\delta^2V\\
\chi_2&=& \overline{D}^2\overline\delta V\\
\chi_3&=& -\overline{D}^3V\\
\chi_4&=& D^2\delta V\\
\chi_5&=& D\delta^2 V
\end{eqnarray*}
The Laplace equation on $V$ and stationarity forces relations among these, since
\[D\overline{D}V+\delta\overline{\delta}V=0=(D+\overline{D})V.\]
We may write the components of $\chi$ in terms of three complex numbers $\lambda, \mu, \nu$ and a real $a$ as
\[\chi_0=\nu;\; \chi_1=-\overline\chi_5=\mu;\; \chi_2=\overline\chi_4=\lambda;\; \chi_3=ia.\]
Also at $p$
\[T^{EE'}\nabla_{EE'}\psi_{ABCD}=0\]
which forces relations between the components of $\chi_{ABCDE}$ and its Hermitian conjugate:
\begin{eqnarray*}
\chi^\dagger_0&=&\chi_1\\
\chi^\dagger_1&=&\chi_2\\
\chi^\dagger_2&=&\chi_3\\
\chi^\dagger_3&=&\chi_4\\
\chi^\dagger_4&=&\chi_5
\end{eqnarray*}
with also
\[
\chi^\dagger_5=\overline\chi_0.
\]
from Hermiticity of  $\nabla_{E'E}\psi_{ABCD}$.

We expand $\omega_{A'B'}$ and $\rho_{AB}$ in the dyad as
\bea
\omega_{A'B'}&=&\omega_2\tilde{o}_{A'}\tilde{o}_{B'}-
2\omega_1\tilde{o}_{(A'}\tilde{o}^{\dagger}_{ B')}+\omega_0\tilde{o}^{\dagger}_{ A'}\tilde{o}^{\dagger}_{ B'}\\
\rho_{AB}&=&\rho_2o_Ao_B-2\rho_1o_{(A}o^\dagger_{B)}+\rho_0o^\dagger_{A}o^\dagger_{B}
\eea
and take dyad components of (\ref{cond2}) to obtain the system
%

\be
\label{sys2}
\left( \begin{array}{ccccc}
 -5\mu & 5\nu & 0 & -\nu & -\mu\\
-4\lambda & 3\mu & \nu & -\mu & -\lambda\\
-3ia & \lambda & 2\mu & -\lambda & -ia\\
-2\overline\lambda & -ia & 3\lambda & -ia & -\overline\lambda\\
\overline\mu & -3\overline\lambda & 4ia & -\overline\lambda &\overline\mu\\
0 & 5\overline\mu & 5\overline\lambda & \overline\mu & -\overline\nu
 \end{array} \right) 
\left(\begin{array}{c}
\rho_0\\
\rho_1\\
\rho_2\\
\omega_1/8\\
\omega_2/8
\end{array} \right) =0.
\ee
This is a system of six (complex) equations in five unknowns, since $\omega_0$ is absent (but recall that reality of $\omega_{A'B'}$ entails $\omega_0=\overline\omega_2$). We may conclude that the only solution of (\ref{sys2}) has zero $\omega_{A'B'}(p)$, if the rank of the coefficient matrix is five. This is generically true, as we see by taking the special case $\nu=\lambda=0$ and $\mu=1$. 
We omit the fourth row of (\ref{sys2}) and calculate the determinant of the resulting $5\times 5$ matrix to be $64(1-a^2)$ . This is nonzero if $a^2\neq 1$, and then the only solution has $\omega_{A'B'}(p)=0$. 
We need to check back to (\ref{co2}), and with these choices we do find that we may deduce $K^i=0$ provided $a^2\neq 1$. Thus there is an open set of choices for $V$ which force  $\omega_{A'B'}(p)=0$ and provide ASD metrics which are not conformal to K\"ahler in any neighbourhood of the chosen point $p$.
\section{The relation with projective invariants}
\label{section_3}
A two--dimensional  projective structure on $(U, [\gamma])$  is an equivalence class of torsion--free
connections on $TU$
\[
\gamma_{A'B'}^{C'}\cong\gamma_{A'B'}^{C'}+\delta_{A'}^{C'} 
\beta_{B'}+\delta_{B'}^{C'}\beta_{A'}, \quad
A', B', C'=1, 2.
\]
(the reasons for using  primed indices will become clear shortly)
where $\beta_{A'}$ is a one--form on $U$.

Consider  the symmetric projective connection $\nabla^{\Pi}$ with connection
symbols
\[
\label{thomas_symbols}
\Pi_{A'B'}^{C'}=\gamma_{A'B'}^{C'}-\frac{1}{3}\gamma_{D'A'}^{D'}
\delta_{B'}^{C'}-\frac{1}{3}\gamma_{D'B'}^{D'}\delta_{A'}^{C'}.
\]
These symbols  do not 
depend on a choice of $\gamma$ in the projective class.
Let $x^{A'}$ be local coordinates on $U$ and let $(x^{A'}, z^{A'})$ be local 
coordinates on $TU$. Define a  $(2, 2)$ conformal structure on $TU$ by
\be
\label{conf_lift}
g=dz_{A'}\otimes dx^{A'}  - \Pi^{C'}_{A'B'}\,z_{C'}\, dx^{A'}\otimes dx^{B'}.
\ee
This is a projectively invariant version of the Riemannian extension
studied by Walker \cite{walker} (also introduced in \cite{YI73}
as a horizontal lift).

This conformal structure  
is anti--self--dual and admits a twisting, null conformal  Killing vector
$K=z^{A'}\p/\p z^{A'}$. Thus it fits into the classification
\cite{DW07} as explained in \cite{BDE08}. Choose a spin frame\footnote{
In this frame 
$K^a=o^{A'}\iota^{A}$ where $\nabla_{AA'}\iota^{B}=0$ and $o^{A'}=z^{A'}$. 
The conformal Killing vector has an ASD derivative $dK=dz_{A'}\wedge dx^{A'}$ 
and
a non-zero twist. It defines two integrable distributions
\[
\{K, {\mathcal S}=z^{A'}\frac{\p}{\p x^{A'}}-\Pi_{A'B'}^{C'}z^{A'}z^{B'}\frac{\p}{\p z^{C'}}\}\qquad \mbox{and}\qquad \{\p/\p z^{1'}, \p/\p z^{2'}\}
\]
where ${\mathcal S}$ is the geodesic spray of the projective structure.}
\[
\nabla_{A'1}=\frac{\p}{\p z^{A'}}, \qquad
\nabla_{A'0}=\frac{\p}{\p x^{A'}}-\Pi_{A'B'}^{C'}z^{B'}\frac{\p}{\p z^{C'}}
\]
and set $\delta_{A'}=\iota^{A}\nabla_{AA'}=\p/\p z^{A'}$. 
Then the spin connection, and curvature
in terms of projective curvature of the projective covariant derivative
$\nabla^{\Pi}$ on $U$ are
\begin{eqnarray}
\Gamma_{AB}&=&\iota_{A}\iota_{B}\Rho_{A'B'}z^{A'} dx^{B'},\quad
\Gamma_{A'B'}=\Pi_{A'B'}^{C'} dx_{C'},\quad \Lambda=0,\\
\phi_{ABA'B'}&=&\iota_{A}\iota_{B}\Rho_{A'B'},\quad
\psi_{A'B'C'D'}=0, 
\quad \psi_{ABCD}=\iota_{A}\iota_{B}\iota_{C}\iota_{D}(z^{A'}Y_{A'})
\end{eqnarray}
where $\Rho_{A'B'}$ is the symmetric Ricci tensor of the 2D projective structure
and 
\[
Y_{C'}=\varepsilon^{A'B'}Y_{A'B'C'}, \qquad  
Y_{A'B'C'}=\frac{1}{2}(\nabla^\Pi_{A'}\Rho_{B'C'}-\nabla^{\Pi}_{B'}\Rho_{A'C'})
\]
is the  Cotton tensor of $[\gamma]$. The components of $Y_{C'}$
are given by (2.14) in \cite{BDE08} and referred to as Liouville's invariants.
In particular $(U, [\gamma])$ is flat if $Y_{C'}=0$. Note
that the spinor $\iota^{A}=(0, 1)$ is covariantly constant on $TU$.

In \cite{BDE08} it was shown that a projective structure
comes from a (possibly Lorentzian) metric on $U$ if 
and only if
there
exists a covariantly constant section $(\sigma^{A'B'}, \mu^{A'}, \rho)$
of  a connection
\begin{equation}
\label{prolonged}\left\lgroup\begin{array}c
\sigma^{B'C'}\\[3pt]
\mu^{B'}\\[3pt]
\rho
\end{array}\right\rgroup\stackrel{{\mathcal{D}}_{A'}}{\longmapsto}
\left\lgroup\begin{array}c
\nabla_{A'}\sigma^{B'C'}-\delta_{A'}^{B'}\mu^{C'}-
\delta_{A'}^{C'}\mu^{B'}\\[3pt]
\nabla_{A'}\mu^{B'}-\delta_{A'}^{B'}\rho+\Rho_{A'C'}
\sigma^{B'C'}\\[3pt]
\nabla_{A'}\rho+2\Rho_{A'B'}\mu^{B'}-2Y_{A'B'C'}\sigma^{B'C'}
\end{array}\right\rgroup,\end{equation}
on a rank $6$ vector bundle over $U$ for which $\sigma^{A'B'}=\sigma^{(A'B')} 
$ is non-degenerate. 
This condition is projectively invariant (when appropriate projective
weights are used) so is also true when the covariant
derivative $\nabla$ with respect to a representative
$\gamma\in [\gamma]$ is replaced by an invariant 
derivative $\nabla^{\Pi}$.
Given such a section, the contravariant metric is constructed by
\be
\label{metric_from_sigma}
h^{A'B'}=\det{(\sigma)}\sigma^{A'B'}.
\ee
The necessary condition for the existence
of the parallel section is obtained by commuting derivatives.
It gives
\be
\label{2Dcond}
5Y_{A'}\,\mu^{A'}+\nabla^{\Pi}_{A'}Y_{B'}\,\sigma^{A'B'}=0.
\ee
(This is (7.46) or (3.20) in \cite{BDE08}).

Now form a 10-tractor $(\omega^{A'B'}, K^{AA'}, \rho^{AB})$ 
from the 6-tractor  $(\sigma^{A'B'}, \mu^{A'}, \rho)$ by
\be
\label{10_from_6}
\omega^{A'B'}=\sigma^{A'B'}, \qquad K^{AA'}=\iota^{A}\mu^{A'}, \qquad
\rho^{AB}=\iota^{A}\iota^{B}\rho.
\ee
Then the first and the last term in the condition (\ref{cond})
vanish. Using
\[
\iota^{A}\nabla_{AA'}(z^{B'}Y_{B'})=Y_{A'}, \qquad
\nabla_{A(A'}Y_{B')}=\iota_{A}\nabla^{\Pi}_{(A'}Y_{B')}
\]
reduces the five conditions (\ref{cond}) to one condition
\[
\iota_{A}\iota_{B}\iota_{C}\iota_{D}(5Y_{A'}\mu^{A'}+\nabla^{\Pi}_{A'}Y_{B'}\,\sigma^{A'B'})=0
\]
which holds if (\ref{2Dcond}) does.
Thus, given the conformal structure (\ref{conf_lift}), the
6-tractor bundle with connection over $U$ embeds in a 10-tractor bundle
with connection over $M=TU$ and the rank $5$ curvature of the latter
is given by a rank $1$ curvature of the former.

Differentiating (\ref{cond}) gives
$25$ conditions on $10$ unknowns so some constraints must hold
for the conformal structure. But we know that
these will hold automatically for (\ref{conf_lift}): in \cite{BDE08}
it was shown that the first constraint arises
after taking two derivatives of (\ref{2Dcond}). There is perhaps no surprise
here  - (\ref{conf_lift}) is type $N$ and the lowest order
obstructions tend to vanish in this case.

A point of caution is needed: if the projective structure
is metrisable by a Riemannian metric, then 
(\ref{10_from_6}) implies that $\omega_{A'B'}\omega^{A'B'}>0$
and thus $\omega_{A'B'}$ gives rise to a K\"ahler metric, albeit
in (2, 2) signature. If on the other hand, the metric
underlying the projective structure is Lorentzian
then $|\omega|^2<0$ and one instead obtains a para--K\"ahler structure: there exists a (2, 2) metric $g$ and 
an almost--product structure 
\[
J:TM\longrightarrow TM,\qquad J^2=\mbox{Id}
\]
such that
\begin{itemize}
\item The structure $J$ is integrable in the sense that the eigenspaces of $TM$ corresponding to eigenvalues $\pm 1$ of $J$ are integrable distributions.
\item $g(JX, JY)=-g(X, Y)$ for all $X, Y\in TM$.
\item The two--from $\Sigma:=g(J\;.,\;.)$ is closed.
\end{itemize}

We have proved the `only if' part of the following
\begin{theo}
\label{theo_projective}
The projective structure $(U, [\gamma])$ is metrisable
if and only if its Riemannian extension {\em(\ref{conf_lift})\em} 
contains a K\"ahler or a para--K\"ahler 
metric in its conformal class.
\end{theo}
{\bf Proof.} It remains to prove the `if' part,
and show that if $(U, [\gamma])$ is metrisable then the $(2, 2)$ metric
(\ref{conf_lift}) on $TU$ is K\"ahler. 
Let $h$ be a metric on $U$. First assume that $h$ is Riemannian. 
Its conformal class defines a complex structure $j:TU\rightarrow TU, 
j^2=-\mbox{Id}$.
Let $\omega $ be a canonical symplectic structure on $T^*U$ and let
\[
T(TU)=V\oplus H
\]
be the splitting of the tangent space to $TU$ 
into vertical and horizontal 
components with respect to the Levi--Civita connection of $h$. The complex
structure $J$ on $TU$ defined   
by taking the complex structure $j$  on each factor $H$ and $V$.

We regard $h$
as an isomorphism between $TU$ and $T^*U$, 
and define a metric $g$
on  $TU$ by
\be
\label{lichn}
g(X, Y)=h^*(\omega)(JX, Y), 
\ee
This agrees with (\ref{conf_lift}) if local coordinates are 
adapted. To see it use the definition of $J$ to find
\[
J\Big(\frac{\p}{\p z^{A'}}\Big)=j^{B'}_{A'}\frac{\p}{\p z^{B'}},\quad
J\Big(\frac{\p}{\p x^{A'}}\Big)=j^{B'}_{A'}\frac{\p}{\p x^{B'}}+
z^{B'}(\gamma_{A'B'}^{C'}j^{D'}_{C'}-\gamma_{C'B'}^{D'}j^{C'}_{A'} )\frac{\p}{\p z^{D'}},
\]
where $\gamma$ is the Levi--Civita connection of $h$.
A triple $(g, J, h^*(\omega))$ is a K\"ahler structure on $TU$.

If the metric $h$ on $U$ is Lorentzian, then its conformal 
structure defines a product structure $j$ on $TU$ with $j^2=\mbox{Id}$. The argument given in the proof still applies, but
it leads to a product structure $J$ on $TM=T(TU)$ and eventually
to a para--K\"ahler structure on $M$.
\koniec
\subsection{Anti--self--dual null K\"ahler structures}
\label{subsection_31}
In this Section we shall invariantly characterise the Riemannian extensions
as a subclass of all (2, 2) ASD conformal structures 
which admit a parallel real section of $\spp$.

In \cite{D02} (2, 2) ASD metrics which admit a
covariantly constant real spinor $\iota^{A'}$ were studied. 
These
were called null K\"ahler structures as the  endomorphism
$N^a_b=\iota^{A'}\iota_{B'}\varepsilon^{A}_B$ satisfies
\be
\label{null_kahler}
N^2=0, \qquad g(NX, Y)+g(X, NY)=0, \qquad \nabla N=0.
\ee
It resembles the K\"ahler condition albeit null.

The condition
(\ref{null_kahler}) is equivalent to the existence of parallel
$\iota^{A'}$ or $\iota^{A}$ and
in the following we shall
choose the spinor to be $\iota^{A}$, so 
that (\ref{null_kahler}) holds with
\[
N^a_b=\iota^{A}\iota_{B}\varepsilon^{A'}_{B'}.
\]
Given a null K\"ahler structure there exist a local coordinate
system $(x^{A'}, z^{A'})$ and a function $\Theta=\Theta(x^{A'}, z^{A'})$
such that 
\[
g=dz_{A'}\otimes dx^{A'}+\frac{\p ^2 \Theta}{\p z^{A'}\p z^{B'}} dx^{A'} \otimes dx^{B'}, \qquad
N=dx^{A'}\otimes\frac{\p}{\p z^{A'}}
\]
where the indices are raised and lowered 
using $\varepsilon_{A'B'}$ \cite{B00, D02}.

The  self--duality conditions 
$\psi_{ABCD}=0$ imposed on $g$ lead to a fourth
order integrable PDE for $\Theta$. This was shown in 
\cite{D02}, where the opposite orientation was used.
If on the other hand the anti-self-duality conditions 
$\psi_{A'B'C'D'}=0$ are imposed 
$\Theta$ can be found explicitly and
we obtain the following
\begin{prop}
There is a one--to--one correspondence between ASD Null K\"ahler structures
where the parallel spinor is a section of $\spp$, and Riemannian extensions
of projective structures of the form {\em(\ref{conf_lift})}.
\end{prop}
{\bf Proof.}
Choose a spin frame
\[
\nabla_{A'1}=\iota^{A}\nabla_{AA'}=\frac{\p}{\p z^{A'}}, \quad
\nabla_{A'0}=o^{A}\nabla_{AA'}=\frac{\p}{\p x^{A'}}
+\frac{\p \Theta}{\p z^{A'}\p z^{B'}}\frac{\p}{\p z_{B'}},
\]
and  set
\[
f=\frac{\p^2\Theta}{\p x^{A'}\p z_{A'}}+\frac{1}{2}
\frac{\p^2 \Theta}{\p z^{A'}\p z^{B'}}\frac{\p ^2 \Theta}{\p z_{A'}\p z_{B'}}.
\]
Then
\begin{eqnarray*}
\psi_{A'B'C'D'}&=&\delta_{A'}\delta_{B'}\delta_{C'}\delta_{D'}\Theta,\qquad
\psi_{ABCD}=\iota_{A}\iota_{B}\iota_{C}\iota_{D}\square f,\qquad
\phi_{ABA'B'}=\iota_{A}\iota_{B}\delta_{A'}\delta_{B'}f,\\
\Gamma_{A'B'}&=&\delta_{A'}\delta_{B'}\delta_{C'}\Theta\, dx^{C'},
\qquad \Gamma_{AB}=\iota_{A}\iota_{B}\delta_{A'}f d x^{A'},\qquad
\Lambda=0,
\end{eqnarray*}
where
\[\delta_{A'}=\frac{\p}{\p z^{A'}}, \qquad\mbox{and}\qquad
\square=\frac{\p^2}{\p x^{A'}\p z_{A'}}+\frac{\p ^2 \Theta}{\p z_{A'}\p z_{B'}}
\frac{\p^2}{\p z^{A'}\p z^{B'}}.
\]
Therefore the self--duality condition implies
that
\[
\Theta=-\frac{1}{6}\Pi_{A'B'}^{C'}z^{A'}z^{B'}z_{C'}
\]
for some $\Pi_{A'B'}^{C'}=\Pi_{A'B'}^C(x^{D'})$ such that
\[
\Pi_{A'B'}^{C'}=\Pi_{B'A'}^{C'}, \qquad \Pi_{A'B'}^{A'}=0
\]
(the terms of order lower
than 3 in $z^{A'}$ can be eliminated by redefining $\Theta$ and translating 
the coordinates $z^{A'}\rightarrow z^{A'}+t^{A'}(x^{B'})$).  
Now $f=(1/2)P_{A'B'} z^{A'}z^{B'}$ and the curvature
and connection coefficients agree with those of 
(\ref{conf_lift}).\koniec
\section{Twistor Theory.} We shall end the paper by briefly describing the
twistor origins of the `conformal to K\"ahler' obstructions. 

Let ${B}$ be a twistor--space (a complex three-fold with
an embedded rational curve with normal bundle 
${\mathcal O}(1)\oplus{\mathcal O}(1)$) corresponding to an ASD conformal structure $(M, [g])$ \cite{Pe76, AHS78}. A K\"ahler structure in $[g]$ corresponds to a preferred section of the
anti-canonical divisor bundle
${\kappa_B}^{-1/2}$, where $\kappa_B$ is the canonical bundle of 
$B$ \cite{pontecorvo}. Restriction of  $\kappa_B$ to a rational curve gives
a line bundle isomorphic to the fourth power of the tautological line bundle
$\mathcal{O}(-1)$. Therefore 
the restriction of the preferred  section to a curve, 
pulled back to the total space of $\spp'\rightarrow M$ is of the form
$\pi^{A'}\pi^{B'}\omega_{A'B'}$, where $\pi^{A'}$ are coordinates on the fibres
of $\spp'$ and $\omega_{A'B'}$ satisfies the conformally invariant linear 
equation (\ref{ck1}).  

In Section \ref{subsection_ASD}
we  constructed a rank--10 vector bundle $E\rightarrow M$ with connection,
such that the parallel sections of this bundle correspond to solutions
to (\ref{ck1}), and therefore to sections of ${\kappa_B}^{-1/2}$.
The forward Ward transform \cite{Wa77,AHS78} of $E$ gives rise to a rank--10 holomorphic vector
bundle $\mathcal{E}\rightarrow B$ (with no connection) which is holomorphically trivial on the twistor curves. This holomorphic vector bundle
can be also constructed directly from the twistor data and is given by
the second-jet bundle $\mathcal{E}=J^2({\kappa_B}^{-1/2})$.


\begin{thebibliography}{99}

\bibitem{AHS78} Atiyah, M.F., Hitchin, N. J. \& Singer, I.M. (1978)
Self-duality in four-dimensional Riemannian geometry. 
Proc. Lon. Math. Soc {\bf A 362}, 425-461.

\bibitem{AG97} Apostolov, V.\and Gauduchon, P. (1997) 
The Riemannian Goldberg-Sachs theorem.  Internat. J. Math.  {\bf 8}, 
421--439.


\bibitem{BEG}
Bailey, T. N., Eastwood, M. G. \and Gover, A. R. (1994) Thomas's structure bundle for conformal, projective and related structures.  Rocky Mountain J. Math.  {\bf{24}}, 1191--1217.
%
\bibitem{B00} Bryant, R. L. (2000)
Pseudo-Riemannian metrics with parallel spinor fields and vanishing
Ricci tensor. Global analysis and harmonic analysis, 
53-94, Semin. Congr., 4, Soc. Math. France, Paris.

\bibitem{BDE08} Bryant, R. L., Dunajski, M. \& Eastwood, M.
(2008) Metrisability of two-dimensional projective structures
{\tt arXiv:0801.0300v1}, to appear in Journal of Differential Geometry.

\bibitem{ChValT}
Chave, T., Valent, G. and Tod, K.P. (1996)
$(4,0)$ and $(4,4)$ sigma models with a tri-holomorphic Killing vector.  
Phys. Lett. {\bf B 383}  262--270. 

\bibitem{Der83} Derdzi\'nski, A. (1983)  
Self-dual Kähler manifolds and Einstein manifolds of dimension four.  Compositio Math.  {\bf 49}, 405--433. 

\bibitem{D02}   Dunajski, M. (2002) 
Anti-self-dual four-manifolds with a parallel real spinor,  
Proc. Roy. Soc. Lond. {\bf A 458}, 1205-1222.

\bibitem{DW07}
Dunajski, M., \& West, S. (2007)
Anti-self-dual conformal structures from projective structures.
{ Comm. Math. Phys. {\bf 272}}, 85--118.


\bibitem{KLP97} Kim, J., LeBrun, C. \& 
Pontecorvo, M. (1997)
Scalar-flat K\"ahler surfaces of all genera. 
J. Reine Angew. Math. {\bf 486} , 69--95.


\bibitem{LeB}
LeBrun, C. R. (1982)
${\mathcal H}$-space with a cosmological constant.
Proc. Roy. Soc. London Ser. {\bf{A 380}} , no. 1778, 171--185. 


\bibitem{Lebrun86} LeBrun, C. (1986)
On the topology of self-dual $4$-manifolds.
Proc. Amer. Math. Soc. {\bf 98}  637--640.


\bibitem{L91} LeBrun, C.R. (1991)
Explicit self-dual metrics on 
$\CP^2 \# \cdots \# \CP^2$, J. Diff. Geom. {\bf 34}  233-253.

\bibitem{LM08}
LeBrun, C. \& Maskit, B. (2008)
On optimal 4-dimensional metrics. 
J. Geom. Anal. {\bf 18}  537--564.


\bibitem{SN}  Nurowski, P., Sparling, G. A. J. 2003
Three-dimensional Cauchy-Riemann structures and second-order
ordinary differential equations, \emph{Class. Quant. Grav.}
\textbf{20}, 4995-5016.

\bibitem{Pap}
Papadopoulos G. (1995) Elliptic monopoles and $(4,0)$-supersymmetric sigma models with torsion.
Phys. Lett. {\bf B 356}  249--255. 

\bibitem{Pe76} Penrose, R. (1976) Nonlinear gravitons and curved
twistor theory, Gen. Rel. Grav. {\bf 7}, 31-52.

\bibitem{PR}
Penrose, R. \& Rindler, W. (1987, 1988) 
{\em Spinors and space-time.  Two-spinor calculus and relativistic fields. }
Cambridge Monographs on Mathematical Physics. Cambridge University Press, 
Cambridge

\bibitem{pontecorvo} Pontecorvo, M. (1992)
On twistor spaces of anti-self-dual hermitian surfaces.
Trans. Am. Math. Soc. {\bf 331}, 653--661.

\bibitem{Prz}
Przanowski, M. (1991) Killing vector fields in self-dual, Euclidean Einstein spaces with $\Lambda\neq 0$.  
J. Math. Phys.{\bf 32}  1004--1010.

\bibitem{RS05}
Rollin, Y.,  Singer, M. (2005)
Non-minimal scalar-flat K\"ahler surfaces and parabolic stability. 
Invent. Math. {\bf 162}, 235--270.

\bibitem{S03} Semmelmann, U. (2003) 
Conformal Killing forms on Riemannian manifolds.  Math. Z.  {\bf 245}  
503--527. 

\bibitem{T97} Tod, K. P. (1995) 
The ${\rm SU}(\infty)$-Toda field equation and special four-dimensional metrics.Geometry and physics (Aarhus, 1995),  307--312.
Lecture Notes in Pure and Appl. Math., 184, Dekker, New York, 1997

\bibitem{walker} Walker, A. G. (1953) Riemann extensions of non-Riemannian
spaces. In {\em Convegno di Geometria Differenziale }. Venice.

\bibitem{Wa77} Ward, R.S. (1977) On self-dual gauge fields, Phys.
Lett. {\bf 61A}, 81-2.

\bibitem{YI73}Yano, K. \& Ishihara, S. (1973) 
Tangent and cotangent bundles: differential geometry. Pure and Applied Mathematics, No. 16. Marcel Dekker, Inc., New York


\end{thebibliography}
\end{document}